\documentstyle[twoside,12pt]{article}
\newtheorem{prop}{Proposition}[section]
\newtheorem{dfn}[prop]{Definition}
\newtheorem{theo}[prop]{Theorem}

\newtheorem{rem}[prop]{Remark}
\newtheorem{coro}[prop]{Corollary}

\newtheorem{exam}[prop]{Example}

\title{MIRROR DUALITY AND STRING-THEORETIC HODGE NUMBERS}

\author{Victor V. Batyrev\thanks{Supported by Deutsche
Forschungsgemeinschaft.} \\
\small  Universit\"at-GHS-Essen, Fachbereich  6,  Mathematik \\
\small  Universit\"atsstra{\ss}e 3, 45141  Essen, Germany  \\
\small e-mail: victor.batyrev@aixrs1.hrz.uni-essen.de \\
and \\
Lev A. Borisov\thanks{Supported by Alfred P. Sloan Doctoral
Dissertation Fellowship} \\
\small Department of Mathematics,  University of Michigan \\
\small Ann Arbor, Michigan 48109-1003, USA \\
\small e-mail: lborisov@math.lsa.umich.edu}

\begin{document}

\date{}

\maketitle

\begin{abstract}
We prove in full generality the  mirror duality conjecture
for  string-theoretic Hodge numbers
of Calabi-Yau complete intersections
in Gorenstein toric Fano varieties.
The proof is based on properties of intersection cohomology.
\end{abstract}

\section{Introduction}

The first author has conjectured
that the polar duality of  reflexive polyhedra induces the
mirror involution for Calabi-Yau hypersurfaces in Gorenstein
toric Fano varieties \cite{bat.dual}. The second author has proposed
 a more general duality which conjecturally
induces  the mirror involution
for Calabi-Yau {\em complete intersections}
in Gorenstein toric Fano varieties \cite{borisov}. The most general form
of the combinatorial duality which includes mirror
constructions of physicists for rigid Calabi-Yau manifolds was formulated
by both authors in \cite{batyrev-borisov1}.

The main purpose of our paper is to show that all proposed
combinatorial dualities agree with the
following   Hodge-theoretic property of mirror symmetry
predicted by physicists:
\bigskip

{\em If two {smooth} $n$-dimensional
Calabi-Yau manifolds $V$ and $W$ form a {mirror pair},
then their Hodge numbers satisfy the relation
\begin{equation}
h^{p,q}(V) = h^{n-p,q}(W), \;\;\;  0 \leq p, q \leq n.
\label{H1}
\end{equation}
 }
\bigskip

A verification of this property becomes rather
non-trivial if we don't make restrictions
on the dimension $n$. The main difficulty is the necessity to work
with singular Calabi-Yau varieties whose singularities in
general don't admit any crepant desingularization.
This difficulty was the motivation for introduction
of so called {\em string-theoretic Hodge numbers} $h^{p,q}_{\rm
st}(V)$ for singular $V$ \cite{batyrev.dais}.
The string-theoretic Hodge numbers
$h^{p,q}_{\rm st}(V)$ coincide with the usual Hodge numbers
$h^{p,q}(V)$ if $V$ is smooth, and with the usual Hodge numbers
of a crepant desingularization $\hat{V}$ of $V$ if  such a desingularization
exists. Therefore the property (\ref{H1}) must be
reformulated as follows:
\bigskip

{\em Let $(V,W)$ be  a mirror pair of singular $n$-dimensional
Calabi-Yau varieties. Then the string-theoretic Hodge numbers
satisfy the duality:
\begin{equation}
h^{p,q}_{\rm st}(V) = h^{n-p,q}_{\rm st}(W), \;\;\;  0 \leq p, q \leq n.
\label{H2}
\end{equation}
}
\bigskip

The string-theoretic Hodge numbers
for Gorenstein algebraic varities with toroi\-dal or quotient singularities
were introduced and studied  in \cite{batyrev.dais}. It was also conjectured
in \cite{batyrev.dais} that the  conbinatorial construction
of mirror pairs of Calabi-Yau complete intersections in Gorenstein toric
Fano varieties satisfiies the duality (\ref{H2}). This conjecture has been
proved in \cite{batyrev.dais}
for mirror pairs of Calabi-Yau hypersurfaces of arbitrary
dimension that
can be obtained by the Greene-Plesser construction
\cite{greene.plesser}. Some other
results  supporting  this conjecture have been obtained
in \cite{bat.dual,batyrev-borisov2,roan}. Additional  evidence
in favor of the conjecture has been received by explicit
computations of instanton sums using generalized hypergeometric functions
\cite{batyrev.straten,HKTY,klemm-theisen,LT}.
\bigskip

The paper is organized as follows:
\bigskip

In Section 2,  we introduce some
polynomials  $B(P; u,v)$ of an Eulerian partially ordered set
$P$ using results of Stanley \cite{stanley}.
It seems that the polynomial $B(P; u,v)$  have independent interest
in  combinatorics. For our purposes, their most important
property is the relation between $B(P; u,v)$ and $B(P^*; u,v)$,
where $P^*$ is the dual to $P$ Eulerian poset (Theorem \ref{duality0}).
\bigskip

In Section 3, we give
an explicit formula for the polynomial $E(Z;u,v)$ which describes
the mixed Hodge structure of an affine hypersurface $Z$ in an algebraic
torus ${\bf T}$ (Theorem \ref{f-formula}).
We remark the following: the explicit formula for
$E(Z; 1,1)$ is due to Bernstein,
Khovanski\^i and Kushnirensko \cite{khov,kush}; the computation
of the polynomial $E(Z; t,1)$ which describes
the Hodge filtration on $H^*_c(Z)$
is due to Danilov and Khovanski\^i \cite{danilov.khov} (see also
\cite{bat.mixed}); the polynomial $E(Z; t,t)$ which describes
the weight  filtration on $H^*_c(Z)$ has been computed by Denef and Loeser
\cite{DL}.
\bigskip

In Section 4, we derive an explicit formula for the polynomial
$E_{\rm st}(V;u,v)$ where $V$ is a Calabi-Yau complete intersection in
a Gorenstein toric Fano variety (Theorem \ref{st.formula}).
The  coefficients of
$E_{\rm st}(V;u,v)$ are equal up to a sign to string-theoretic Hodge numbers
of $V$. Since our formula is written in terms of $B$-polynomials
as a sum over pairs of lattice points contained
in the corresponding pair of dual to each other reflexive Gorenstein
cones $C$ and $C^*$, the mirror duality for string-theoretic Hodge
numbers becomes immediate consequence of the duality for $B$-polynomials after
the transposition $C \leftrightarrow C^*$ (Theorem
\ref{duality1}). Following some recent development
of ideas of Witten \cite{witten} by Morisson and Plesser
\cite{morrison.plesser},
we conjecture that the  formula obtained in this paper
gives the spectrum of
the  abelian gauge theory in two dimensions
which could  be constructed from any pair $(C,C^*)$
of two dual to each other reflexive Gorenstein cones.
\bigskip

{\bf Acknowledgements.} The first author would like to thank
for hospitality the University of Warwick where the paper
was finished.
\bigskip

\section{Combinatorial polynomials of Eulerian posets}

Let $P$ be a finite poset (i.e., finite partially ordered set).
Recall that the M\"obius function $\mu_P(x,y)$ of a poset $P$ is a
unique integer valued function on $P \times P$ such that for
every function $f: P \rightarrow A$ with values in an
abelian group $A$ the following {\em M\"obius inversion formula}
holds:
$$
f(y) = \sum_{x \leq y} \mu_P(x,y) g(x), \;\;\;
\mbox{\rm where $\;g(y) = \sum_{x \leq y} f(x)$}.
$$

{}From now on  we always assume that the poset $P$ has a unique minimal
element $\hat{0}$, a unique maximal element $\hat{1}$, and that
every maximal chain of $P$ has the same length $d$ which will
be called the {\em rank of } $P$. For any $x \leq y$ in $P$,
define the interval
\[ \lbrack x, y \rbrack = \{ z \in P\,:\, x
  \leq z \leq y \}. \]
In particular, we have $P = \lbrack \hat{0}, \hat{1}
\rbrack$.
Define the
rank function $\rho\,:\, P \rightarrow \{0,1,\ldots,d \}$ of $P$ by
setting
$\rho(x)$ equal to the length of any saturated chain in the
interval $\lbrack \hat{0}, x \rbrack$.

\begin{dfn}
{\rm \cite{stanley} A poset $P$
as above is said to be {\em Eulerian} if
for any $x \leq y$ $(x,y \in P)$ we have
$$
\mu_P(x,y) = (-1)^{\rho(y) - \rho(x)}.
$$}
\end{dfn}

\begin{rem}
{\rm It is easy to see that any interval $\lbrack x, y \rbrack
\subset P$ in an Eulerian poset $P$ is again an Eulerian poset
with the rank function $\rho(z) - \rho(x)$ for any
$z \in \lbrack x, y \rbrack$. If an Eulerian poset $P$ has rank
$d$, then the dual poset $P^*$ is again an Eulerian poset with
the rank function $\rho^*(x) = d - \rho(x)$. }
\end{rem}

\begin{exam}
{\rm Let $C$ be an $d$-dimensional finite convex polyhedral cone in
${\bf R}^d$ such that $-C \cap C = \{0\} \in {\bf R}^d$. Then the poset
$P$ of faces of $C$ satisfies all above assumptions
with the maximal element  $C$, the  minimal element $\{0\}$, and
the rank function $\rho$ which is equal to the dimension
of the corresponding face.
It is easy to show that $P$ is an Eulerian poset of rank $d$.}
\end{exam}

\begin{dfn}
{\rm \cite{stanley} Let $P = \lbrack \hat{0}, \hat{1} \rbrack$
be an Eulerian poset of rank $d$. Define two polynomials
$G(P,t)$, $H(P,t) \in {\bf Z} [ t]$ by  the following recursive rules:
$$
G(P,t) = H(P,t) = 1\;\; \mbox{\rm if $d =0$};
$$
$$
H(P,t) = \sum_{ \hat{0} <  x \leq  \hat{1}} (t-1)^{\rho(x)-1}
G(\lbrack x,\hat{1}\rbrack, t)\;\; (d>0),
$$
$$
G(P,t) =
 \tau_{ < d/2 } \left(
(1-t)H(P,t) \right) \;\;( d>0),
$$
where $\tau_{ < r }$ denotes the truncation operator ${\bf
Z}\lbrack t \rbrack  \rightarrow {\bf
Z}\lbrack t \rbrack$ which is defined by
\[ \tau_{< r} \left( \sum_i a_it^i \right) = \sum_{i < r}
a_it^i. \]}
\label{GH}
\end{dfn}

\begin{theo} {\rm \cite{stanley} }
Let $P$ be an Eulerian poset of rank $d \geq 1$. Then
$$
H(P,t) = t^{d-1} H(P,t^{-1}).
$$
\label{symm}
\end{theo}

\begin{prop}
  Let $P$ be an Eulerian poset of rank $d \geq 0$. Then
$$
t^d G(P,t^{-1}) = \sum_{\hat{0} \leq x \leq \hat{1}}
(t-1)^{\rho(x)} G(\lbrack x,\hat{1} \rbrack, t).
$$
\label{inv}
\end{prop}

\noindent
{\em Proof.} The case $d=0$ is obvious. Using \ref{symm},
we obtain
$$
(t-1) H(P,t) = t^d G(P,t^{-1}) -  G(P,t)\;\;(d >0).
$$
Now the statement follows from the formula for $H(P,t)$ in
\ref{GH}. \hfill $\Box$

\begin{dfn}
{\rm  Let $P$ be an Eulerian poset of rank $d$. Define
the polynomial $B(P; u,v) \in {\bf Z}[ u,v]$
by  the following recursive rules:
$$
B(P; u,v) = 1\;\; \mbox{\rm if $d =0$},
$$
$$
\sum_{\hat{0} \leq x \leq \hat{1}}
B(\lbrack \hat{0}, x \rbrack; u,v) u^{d  - \rho(x)}
G(\lbrack x ,  \hat{1}\rbrack, u^{-1}v) = G(P ,uv).
$$ }
\label{Q}
\end{dfn}

\begin{exam}
{\rm Let  $P$ be the boolean algebra of rank $d \geq 1$. Then
$G(P,t) = 1$, $H(P, t) = 1 + t + \cdots + t^{d-1}$, and
$B(P; u,v) = (1-u)^{d}$. }
\end{exam}

\begin{exam}
{\rm Let $C \subset {\bf R}^3$
 be a $3$-dimensional finite convex polyhedral cone with  $k$ $1$-dimensional
 faces ($-C \cap C = \{0\} \in {\bf R}^3$),
 $P$ the Eulerian poset of faces of $C$.
Then $G(P,t) = 1 + (k-3)t$, $H(P, t) = 1 + (k-2)t + (k-2)t^2 +
t^3$, and
$$B(P; u,v) = 1  - (k- (k-3)v)u + (k -(k-3)v)u^2 - u^3.$$
We notice that $B(P; u,v)$ satisfies the
relation
$$B(P; u,v) =(-u)^3B(P; u^{-1},v)$$
which is a consequence of the selfduality $P \cong P^*$ and
a more general property \ref{duality0}.}
\end{exam}

\begin{prop}
Let $P$ be an Eulerian poset of rank  $d> 0$. Then
$B(P; u,v)$  has the following properties:

{\rm (i)} $B(P; u,1) = (1-u)^d$ and  $B(P; 1,v) = 0$;

{\rm (ii)} the degree of
$B(P; u,v)$
with respect to $v$ is less than $d/2$.
\label{degree}
\end{prop}

\noindent
{\em Proof.} The statement (i) follows immediately from \ref{inv} and
the recursive definition of $B(P; u,v)$.  In order to prove
(ii) we use  induction on $d$. By assumption, the degree
of $B(\lbrack \hat{0}, x \rbrack; u,v)$ with respect to $v$ is less than
$\rho(x)/2$. On the other hand, the $v$-degree of
$G(\lbrack x ,  \hat{1}\rbrack; u^{-1}v)$  is less
than $(d- \rho(x))/2$ (see \ref{GH}). It remains to apply the
recursive formula of  \ref{Q}.  \hfill $\Box$

\begin{prop}
Let $P$ be an Eulerian poset of rank $d$. Then  $B$-polynomials of
intervals $[\hat{0},x ]$ and $[x, \hat{1}]$ satisfy the following
relation:
$$
\sum_{ \hat{0} \leq x \leq \hat{1}}
B([\hat{0},x ];u^{-1},v^{-1})(uv)^{\rho(x)}(v-u)^{d-\rho(x)}
=
\sum_{ \hat{0} \leq x \leq \hat{1}}
B([x, \hat{1}];u,v)(uv-1)^{\rho(x)}.$$
\label{relation}
\end{prop}

\noindent
{\em Proof}. Let us substitute $u^{-1},v^{-1}$ instead of $u,v$
in the recursive relation \ref{Q}. So
we obtain
\begin{equation}
\sum_{\hat{0} \leq x \leq \hat{1}} B([\hat{0},x];u^{-1},v^{-1})
u^{-d + \rho(x)}G([x, \hat{1}],
uv^{-1}) = G(P, u^{-1}v^{-1}).
\label{eq1}
\end{equation}
By \ref{inv}, we have
\begin{equation}
G(P,u^{-1}v^{-1}) =  (uv)^{-d} \sum_{\hat{0} \leq x \leq \hat{1}}
(uv-1)^{\rho(x)} G(\lbrack x,\hat{1} \rbrack, uv)
\label{eq2}
\end{equation}
and
$$G([x, \hat{1}],uv^{-1}) =  \sum_{x \leq y \leq \hat{1}}
(u^{-1}v -1)^{\rho(y)-\rho(x)}u^{d-\rho(x)}v^{\rho(x) -d}
G(\lbrack y,\hat{1} \rbrack, u^{-1}v)$$
\begin{equation}
=   \sum_{x \leq y \leq \hat{1}}
u^{d-\rho(y)}v^{\rho(x) -d}(v -u)^{\rho(y)-\rho(x)}
G(\lbrack y,\hat{1} \rbrack, u^{-1}v).
\label{eq3}
\end{equation}
By \ref{Q}, we also have
\begin{equation}
G([x, \hat{1}],uv) = \sum_{x \leq y \leq \hat{1}} u^{d - \rho(y)}
B([x,y];u,v) G([y,\hat{1}],u^{-1}v).
\label{eq4}
\end{equation}
By substitution (\ref{eq4}) in (\ref{eq2}), and two equations
(\ref{eq2}), (\ref{eq3}) in (\ref{eq1})
we obtain:
$$
\sum_{ \hat{0} \leq x \leq y \leq \hat{1}}
B([\hat{0},x ];u^{-1},v^{-1})u^{\rho(x) - \rho(y)}
v^{\rho(x) -d}(v-u)^{\rho(y)-\rho(x)}
G([y, \hat{1}], u^{-1}v) $$
\begin{equation}
 = \sum_{ \hat{0} \leq x \leq y \leq \hat{1}}
B([x, y];u,v)u^{-\rho(y)}v^{-d}(uv-1)^{\rho(x)}G([y, \hat{1}], u^{-1}v).
\label{eq5}
\end{equation}

Now we use  induction on $d$. It is easy to see that
the  equation (\ref{eq5}) and the
induction hypothesis for $y < \hat{1}$
immediately imply  the statement of the proposition. \hfill $\Box$

\begin{prop}
The $B$-polynomials are uniquely determined by
the relation \ref{relation}, by the property of $v$-degree from
\ref{degree}(ii),
and by the initial condition $B(P;u,v) =1$ if $d =0$.
\label{uniq}
\end{prop}

\noindent
{\em Proof.}
Indeed, if we know $B([x,y];u,v)$ for all $\rho(y) - \rho(x) < d$,
then we know all terms in \ref{relation} except for $B(P;u,v)$
on the right hand side and $B(P;u^{-1},v^{-1})
(uv)^{d}$ on the left hand side. Because the $v$-degree
of $B(P;u,v)$ is less than $d/2$, the
possible degrees of monomials with respect to  variable $v$ in
$B(P;u,v)$ and $B(P;u^{-1},v^{-1})(uv)^d$ do not coincide. This
allows us to determine $B(P;u,v)$ uniquely. \hfill
$\Box$

\begin{theo}
Let $P$ be an Eulerian poset of rank $d$,  $P^*$ be the dual
Eulerian poset. Then
$$
B(P; u,v) = (-u)^d B(P^*;u^{-1},v).
$$
\label{duality0}
\end{theo}

\noindent
{\em Proof.} We set
$$
{Q}(P; u,v) = (-u)^d B(P^*;u^{-1},v).
$$
It is clear that ${Q}(P; u,v) =1$ and $v$-degree of
${Q}(P; u,v)$ is the same as $v$-degree of ${B}(P; u,v)$.
By \ref{uniq}, it remains to establish the same
recursive relations for ${Q}(P; u,v)$ as for
${B}(P; u,v)$ in \ref{relation}. The last property follows
from straightforward computations. Indeed, the  equality
\begin{equation}
\sum_{ \hat{0} \leq x \leq \hat{1}}
{Q}([\hat{0},x ];u^{-1},v^{-1})(uv)^{\rho(x)}
(v-u)^{d-\rho(x)} = \sum_{ \hat{0}\leq x \leq
\hat{1}} {Q}([x, \hat{1}];u,v)(uv-1)^{\rho(x)}
\label{rel-dual}
\end{equation}
is equivalent to the relation \ref{relation} for $B(P^*;u,v^{-1})$:
$$
\sum_{ \hat{0}\leq x \leq
\hat{1}}
{B}([x,\hat{1} ]^*;u^{-1},v)(uv^{-1})^{d-\rho(x)}
(v^{-1}-u)^{\rho(x)} $$
$$= \sum_{ \hat{0}\leq x \leq
\hat{1}}
{B}([\hat{0},x]^*;u,v^{-1})(uv^{-1}-1)^{d-\rho(x)},
$$
because
$${Q}([x,\hat{1} ];u,v)  =(-u)^{d-\rho(x)}
 {B}([x,\hat{1} ]^*;u^{-1},v)
$$
and
$${Q}([\hat{0},x];u^{-1},v^{-1})  =(-u)^{-\rho(x)}
 {B}([\hat{0},x]^*;u,v^{-1}).
$$
\hfill $\Box$
\bigskip

\section{E-polynomials of toric hypersurfaces}

Let $M$ and $N$ be two free abelian groups of rank $d$ which are
dual to each other; i.e., $N = {\rm Hom}(M, {\bf Z})$. We denote by
\[ \langle *, * \rangle \;:\; M \times N \rightarrow {\bf Z} \]
the canonical bilinear pairing, and by $M_{\bf R}$ (resp. by $N_{\bf R}$)
the  real scalar extensions of $M$ (resp. of $N$).

\begin{dfn}
{\rm A subset $C \subset M$ is called a
$d$-dimensional rational convex polyhedral cone with vertex
$\{0\} \in M$ if there exists a finite set $\{e_1, \ldots, e_k\}
 \subset M$ such that
\[ C = \{ \lambda_1 e_1 + \cdots + \lambda_k e_k \in M_{\bf R} \; : \;
\mbox{\rm where }\; \; \lambda_i \in {\bf R}_{\geq 0}\; ( i =1,\ldots,k)\}  \]
and $-C + C = M_{\bf R}$, $-C \cap C = \{0\} \in M$. }
\label{cone}
\end{dfn}

\begin{rem}
{\rm If $C \subset M$ is a
$d$-dimensional rational convex polyhedral cone with vertex
 $\{0\} \in M$, then the dual cone
\[ C^* = \{ z \in N_{\bf R} \; : \; \mbox{\rm $\langle e_i, z
\rangle \geq 0$ for all $i \in \{ 1, \ldots, k\}$ } \} \]
is also a $d$-dimensional rational convex polyhedral cone with vertex
$\{0\}$ in the dual space $N_{\bf R}$. Moreover,
there exists a canonical bijective correspondence $F
\leftrightarrow F^*$ between
faces $F \subset C$ and faces $F^* \subset C^*$ $({\rm dim}\, F +
{\rm dim}\, F^* = d)$ :
\[ F \mapsto F^* = \{ z \in C^*\; : \; \mbox{\rm $\langle z', z
\rangle = 0$ for all $z ' \in F$} \} \]
which  reverses  inclusion  relation between
faces. }
\end{rem}

Let $P$ be the Eulerian poset of faces of a $d$-dimensional
rational convex polyhedral cone $C \subset M_{\bf R}$ with
vertex in $\{0\}$. For convenience of notations, we use elements
$x \in P$ as indices and denote
by  $C_x$ the face of $C$ corresponding to $x \in P$, in
particular, we have $C_{\hat{0}} = \{0\}$,
$C_{\hat{1}} = C$, and $\rho(x) = {\rm dim}\, C_x$.
The dual Eulerian poset $P^*$ can be identified with
the poset of faces $C_x^*$ of the dual cone $C^* \subset N_{\bf
R}$.

\begin{dfn}
{\rm A  $d$-dimensional cone $C$ $(d \geq 1)$ as in
\ref{cone}
 is called {\em Gorenstein} if there exists an element
$n_C \in N$ such that
$\langle  z, n_C \rangle >$ for any nonzero $z \in C$ and
all vertices of the $(d-1)$-dimensional convex polyhedron
\[ \Delta(C) = \{ z \in C \; : \; \langle z, n_C \rangle =1\} \]
belong to $M$. This polyhedron
will be called the {\em supporting polyhedron of } $C$.
For convinience, we consider $\{0\}$ as  a
$0$-dimensional Gorenstein cone
with the supporting polyhedron $\Delta(\{0\}): = \emptyset$.
For any $m \in C \cap M$, we define the {\em degree of} $m$
as
$$ {\rm deg}\,m = \langle m, n_C \rangle.$$ }
\end{dfn}

\begin{rem}
{\rm It is clear that any face $C_x$ of a Gorenstein cone is
again a Gorenstein cone with the supporting
polyhedron
\[ \Delta(C_x) = \{ z \in C_x \; : \; \langle z, n_C \rangle =1\}.
\]}
\end{rem}
\bigskip

Now we recall  standard facts from theory of toric varieties
\cite{danilov,fulton,oda} and fix our notations:

Let ${\bf P}(C)$ be the $(d-1)$-dimensional projective
toric variety associated with a Gorenstein cone $C$. By definition,
$$
{\bf P}(C) = {\rm Proj}\, {\bf C} [ C \cap M ]
$$
where ${\bf C} [ C \cap M ]$ is a graded semigroup algebra over ${\bf C}$
of lattice points $m \in C \cap M$. Each face $C_x \subset C$ of
positive dimension defines
an irreducible projective toric subvariety
$$
{\bf P}(C_x) = {\rm Proj}\, {\bf C} [ C_x \cap M ] \subset {\bf P}(C)
$$
which is a compactification of a $(\rho(x) -1)$-dimensional
algebraic torus
$$
{\bf T}_x : = {\rm Spec}\,{\bf C} [ M_x ],
$$
where $M_x \subset M$ is the subgroup of all lattice points $m \in
(-C_x + C_x) \cap M$ such that $\langle m,q \rangle = 0$. Moreover,
the multiplicative group low on ${\bf T}_x$ extends to a regular
action of ${\bf T}_x$ on ${\bf P}(C_x)$ so that
one has the natural stratification
$$
{\bf P}(C_x) = \bigcup_{ \hat{0} < y \leq x} {\bf T}_y
$$
by ${\bf T}_x$-orbits ${\bf T}_y$. We denote by ${\cal O}_{\bf P}(C)(1)$ the
ample tautological sheaf on ${\bf P}(C)$. In particular, lattice points
in $\Delta(C)$ can be identified with a torus invariant basis of
the space of global sections of ${\cal O}_{\bf P}(C)(1)$. We denote by
$\overline{Z}$ the set of zeros of a generic global section of
${\cal O}_{\bf P}(C)(1)$   and set
$$
Z_x := \overline{Z} \cap {\bf T}_x\;\; (\hat{0} < x \leq \hat{1}).
$$
Thus  we have the natural stratification:
$$
\overline{Z} = \bigcup_{ \hat{0} < x \leq \hat{1}} Z_x,
$$
where each $Z_x$  is a smooth affine hypersurface
in ${\bf T}_x$ defined by a generic Laurent polynomial with the Newton
polyhedron $\Delta(C_x)$.

\begin{dfn} Define two functions
$$
S(C_x,t):= (1-t)^{\rho(x)} \sum_{m \in C_x \cap M} t^{{\rm deg}\,m}
$$
and
$$
T(C_x,t):= (1-t)^{\rho(x)} \sum_{m \in Int(C_x) \cap M} t^{{\rm deg}\,m},
$$
where $Int(C_x)$ denotes the relative interior of $C_x \subset C$.
\end{dfn}

The following statement is a consequence of the Serre duality
 (see \cite{danilov.khov,bat.mixed}):

\begin{prop}
$S(C_x,t)$ and $T(C_x,t)$ are polynomials satisfying the relation
$$S(C_x,t) = t^d T(C_x,t^{-1}).$$
\label{ST}
\end{prop}

\begin{dfn}
{\rm \cite{danilov.khov}
Let $X$ be a quasi-projective algebraic variety over ${\bf C}$.
For each pair of integers $(p,q)$, one defines the following
generalization of Euler characteristic:
$$ e^{p,q}(X) = \sum_{k} (-1)^{k} h^{p,q}(H_c^k(X)), $$
where $h^{p,q}(H_c^k(X))$ is the dimension of the $(p,q)$-component
of the mixed Hodge structure of $H_c^k(X)$ \cite{deligne}. The sum
$$E(X; u,v) := \sum_{p,q} e^{p,q}(X) u^p v^q $$
is called {\em $E$-polynomial of $X$}. }
\end{dfn}

Next statement is also due to
Danilov and Khovanski\^i (see \cite{danilov.khov}
\S 4 , or another approach in \cite{bat.mixed}):

\begin{prop}
 We  set $E(Z_{\hat{0}};t,1): = (t-1)^{-1}$.
Then
$$ E(Z_x;t,1) = \frac{(t - 1)^{\rho(x)-1} + (-1)^{\rho(x)}
S(C_x,t)}{t}$$
for $\rho(x) \geq 0$.
\label{s-poly}
\end{prop}

The purpose of this section is to give an explicit formula for
$E$-polynomials of affine hypersurfaces $Z_x \subset {\bf T}_x$.
Following the method of Denef and Loeser \cite{DL}
combined with ideas of Danilov
and Khovanski\^i \cite{danilov.khov}, we compute  $E(Z_x;u,v)$ using
intersection cohomology introduced by Goresky and MacPherson \cite{GM0}.
Recall that intersection cohomology $IH^*(X)$ of a
quasiprojective algebraic variety $X$
of pure dimension $n$ over ${\bf C}$ can be defined as hypercohomology
of the so called {\em intersection complex} $IC^{\bullet}_X$. Moreover,
the intersection complex carries a natural
mixed Hodge structure. The
weight filtration on the $l$-adic version of intersection cohomology
has been introduced and studied  by Beilinson, Bernstein,
Deligne  and Gabber  using theory of perverse sheaves  \cite{BBD}.
The Hodge filtration on intersection cohomology of algebraic varieties
over ${\bf C}$ has been introduced
by M. Saito using his theory of mixed Hodge modules \cite{saito}
(see also \cite{durffe}).

\begin{dfn}
{\rm Let $X$ be a quasi-projective algebraic variety over ${\bf C}$.
For each pair of integers $(p,q)$, one defines the following
generalization of Euler characteristic for intersection cohomology:
$$ e^{p,q}_{\rm int}(X) = \sum_{k} (-1)^{k} h^{p,q}(IH_c^k(X)), $$
where $h^{p,q}(H_c^k(X))$ is the dimension of the $(p,q)$-component
in the mixed Hodge structure of $IH_c^k(X)$. The sum
$$E_{\rm int}(X; u,v) := \sum_{p,q} e^{p,q}_{\rm int}(X) u^p v^q $$
is called {\em intersection cohomology $E$-polynomial of $X$}. }
\end{dfn}

The following statement has been discovered by Bernstein, Khovanski\^i and
Mac\-Pherson (two independent proofs are contained in \cite{DL} and
\cite{fieseler}):

\begin{theo}
$$E_{\rm int}({\bf P}(C);u,v)= H(P,uv) =\sum_{\hat{0} < x \leq \hat{1}}
(uv-1)^{\rho(x)-1}G([x,\hat{1}],uv).$$
Moreover, the cohomology sheaves ${\cal H}^i(IC_{\bf P(C)}^{\bullet})$
are constant on torus orbits ${\bf T}_x$ and
$G([x,\hat{1}],uv)$ is the Poincar{\'e} polynomial describing their
dimensions.
\label{intersection}
\end{theo}

\begin{coro}
Let $\overline{W} \subset {\bf P}(C)$ be a hypersurface that  meets
transversally all toric strata ${\bf T}_x \subset {\bf P}(C)$
that it intersects $($$
\overline{W}$
is not necessary ample$)$. Then
$$E_{\rm int}(\overline{W};u,v)=\sum_{\hat{0} < x \leq \hat{1}}
E(W_x;u,v)G([x,\hat{1}], uv),$$
where $W_x = \overline{W} \cap {\bf T}_x$ $(\hat{0} < x \leq \hat{1})$.
\label{e-form}
\end{coro}

\noindent
{\em Proof.} The statement is essentially cointained in
\cite{DL}(Lemma 7.7). The key fact is that singularities of $\overline{W}$
along $W_x$ are {\em toroidal } (see  \cite{danilov}), i.e. ,
locally isomorphic to toric singularities which appear on ${\bf P}(C)$.
\hfill $\Box$
\bigskip

\noindent
Applying \ref{s-poly}, we obtain:

\begin{coro}
$$E_{\rm int}(\overline{Z};t,1)=\sum_{\hat{0} < x \leq \hat{1}}
\left( \frac{(t-1)^{\rho(x)-1}+(-1)^{\rho(x)}S(C_x,t)}{t} \right)
G([x, \hat{1}], t).$$
\label{s-poly1}
\end{coro}

\begin{dfn}
{\rm Define  $H_{\rm Lef}(P,t)$ to be the polynomial
of degree $(d-2)$ with the following
properties:

(i) $H_{\rm Lef}(P,t) = t^{d-2}H_{\rm Lef}(P,t^{-1})$;

{(ii)} $ \tau_{\leq (d-2)/2} H_{\rm Lef}(P,t) =  \tau_{\leq (d-2)/2}
H(P,t)$.  }
\label{def-lef}
\end{dfn}

\begin{prop}
$$H_{\rm Lef}(P,t) = (1-t)^{-1}
(G(P,t)-t^{d-1}G(P,t^{-1})).$$
\label{p-lef}
\end{prop}

\noindent
{\em Proof.}  Let us set
$$Q(P,t) := (1-t)^{-1}
(G(P,t)-t^{d-1}G(P,t^{-1})).$$
We check that the properties \ref{def-lef}(i)-(ii) are satisfied for $Q(P,t)$.
Indeed \ref{def-lef}(i) follows immediately from the definiton
of $Q(P,t)$. If
$$H (P, t) = \sum_{0 \leq i \leq d-1} h_i t^i$$
and
$$G (P, t) = h_0 +  \sum_{1 \leq i < d/2}(h_i - h_{i-1})t^i,$$
then
$$Q(P,t) = h_0 \frac{1 - t^{d-1}}{1-t} + \sum_{1 \leq i < d/2}
(h_i- h_{i-1}) \frac{t^i - t^{d-1-i}}{1-t}.$$
This shows (ii) and the fact that $Q(P,t)$ is a polynomial.
\hfill $\Box$

\begin{prop}
Define $E_{\rm int}^{\rm prim}(\overline{Z};u,v)$ to be the polynomial
$$
E_{\rm int}^{\rm prim}(\overline{Z};u,v):=
E_{\rm int}(\overline{Z};u,v) - H_{\rm Lef}(P,uv).
$$
Then $E_{\rm int}^{\rm prim}(\overline{Z};u,v)$
is a homogeneous polynomial of degree $(d-2)$.
\label{sum}
\end{prop}

\noindent
{\em Proof.} By the Lefschetz theorem for intersection cohomology
\cite{GM}, we have
isomorphisms
$$ IH^i({\bf P}(C)) \cong IH^i(\overline{Z}), \; \; (0 \leq i < d-2)$$
and the short exact sequence
$$  0 \rightarrow IH^{d-2}({\bf P}(C)) \rightarrow
IH^{d-2}(\overline{Z}) \rightarrow IH^{d-2}_{\rm prim}(\overline{Z})
\rightarrow 0,$$
where $IH^{d-2}_{\rm prim}(\overline{Z})$ denotes the primitive part
of intersection cohomology of $\overline{Z}$ in degree $(d-2)$. By purity
theorem for intersection cohomology \cite{gabber} (see also \cite{durffe}),
the Hodge structure of $IH^{d-2}_{\rm prim}(\overline{Z})$ is pure. On the
other hand, it follows from the Poincar{\'e} duality for
intersection cohomology that
$E_{\rm int}^{\rm prim}(\overline{Z};u,v)$ is the $E$-polynomial of this Hodge
structure.
\hfill $\Box$

\begin{theo}
We set $E(Z_{\hat{0}};u,v):=(uv-1)^{-1}$.
Then $E$-polynomials $E(Z_x;u,v)$
of affine toric hypersurfaces satisfy the
following resursive
relation
$$\sum_{\hat{0} \leq x \leq \hat{1}}
(E(Z_x;u,v)-(uv)^{-1}(uv-1)^{\rho(x)-1})G([x, \hat{1}],uv)$$
$$=v^{d-2} \sum_{\hat{0}  \leq x \leq \hat{1}}
(u^{-1}v) (-1)^{\rho(x)}
S(C_x,uv^{-1}))G([x, \hat{1}],uv^{-1}).$$
\label{recur}
\end{theo}

\noindent
{\em Proof.}
By \ref{s-poly1} and  \ref{p-lef}, we have

\bigskip
\noindent
$ E_{\rm int}^{\rm prim}(\overline{Z};t,1) = E_{\rm int}(\overline{Z};t,1) -
H_{\rm Lef}(P,t) = $
$$= \sum_{\hat{0} < x \leq \hat{1}}
t^{-1}((t-1)^{\rho(x)-1}+(-1)^{\rho(x)}S(C_x,t))
G([x, \hat{1}], t)$$
$$  - (1-t)^{-1}
(G(P,t)-t^{d-1}G(P,t^{-1})).$$
Using \ref{inv}, we obtain
$$ \sum_{\hat{0} < x \leq \hat{1}}
t^{-1}(t-1)^{\rho(x)-1}G([x, \hat{1}], t) =
t^{-1}(t-1)^{-1} (t^d G(P,t^{-1}) - G(P,t)). $$
This yields
\begin{equation}
E_{\rm int}^{\rm prim}(\overline{Z};t,1) =
 \sum_{\hat{0} \leq x \leq \hat{1}}t^{-1}(-1)^{\rho(x)}
S(C_x,t))G([x, \hat{1}],t).
\label{mid1}
\end{equation}

On the other hand, by \ref{e-form} and \ref{p-lef}, we have
\bigskip

\noindent
$E_{\rm int}^{\rm prim}(\overline{Z};u,v) = E_{\rm int}(\overline{Z};u,v)  -
H_{\rm Lef}(P,uv) $ \\
$$= \sum_{\hat{0} < x \leq \hat{1}}
E(Z_x;u,v)G([x, \hat{1}],uv) - (1-uv)^{-1}
(G(P,uv)-(uv)^{d-1}G(P,(uv)^{-1})).$$
Using \ref{inv}, we obtain
$$
\sum_{\hat{0} \leq x \leq \hat{1}}
(uv)^{-1}(uv-1)^{\rho(x)-1}G([x, \hat{1}],uv) =
(uv)^{d-1}(uv-1)^{-1}G(P,(uv)^{-1}).$$
This yields
\begin{equation}
E_{\rm int}^{\rm prim}(\overline{Z};u,v)
= \sum_{\hat{0} \leq x \leq \hat{1}}
(E(Z_x;u,v)-(uv)^{-1}(uv-1)^{\rho(x)-1})G([x, \hat{1}],uv).
\label{mid2}
\end{equation}

\noindent
By \ref{sum},  we have
$$E_{\rm int}^{\rm prim}(\overline{Z};u,v) = v^{d-2}
E_{\rm int}^{\rm prim}(\overline{Z};uv^{-1},1).$$
It remains to  combine (\ref{mid1}) and (\ref{mid2}). \hfill
$\Box$

\begin{dfn}
{\rm Let $m$ be a lattice point in $C \cap M$. We denote
by ${x(m)}$ the minimal element among $x \in P$ such that the face
$C_x \subset C$ contains $m$. The interval  $[ x(m), \hat{1}]
\subset P$ parametrizes  the set of all faces of $C$ containing $m$.
We identify the dual interval $[ x(m), \hat{1}]^*$
with the Eulerian poset  of all
faces $C_x^* \subset C$ such that $\langle m, z \rangle= 0$ for
all $z \in C_x^*$.}
\end{dfn}

\begin{theo}
Let us set $Z: = Z_{\hat{1}}$. Then
there exists the following explicit formula for  $E(Z;u,v)$
in terms of $B$-polynomials:
$$
E(Z;u,v) = \frac{(uv-1)^{d-1}}{uv} + \frac{(-1)^d}{uv}
\sum_{m \in C\cap M}
(v-u)^{\rho(x(m))}B([x(m),\hat{1}]^*; u,v)
\left(\frac{u}{v}\right)^{{\rm deg}\,m}.
$$
\label{f-formula}
\end{theo}

\noindent
{\em Proof.} By induction, $E$-polynomials are uniquely determined
from the recursive formula \ref{recur}. Therefore, it suffices
to show that the functions
$$
\frac{(uv-1)^{\rho(x)-1}}{uv} + \frac{(-1)^{\rho(x)}}{uv}
\sum_{m \in C_x \cap M}
(v-u)^{\rho(x(m))}B([x(m),x]^*; u,v) \left(\frac{u}{v}\right)^{{\rm deg}\,m}
$$
satisfy the same resursive formula as polynomials $E(Z_x; u,v)$.
Indeed, let us  substitute
these functions instead of $E$-polynomials in the left hand side
of \ref{recur} and expand
$$(-1)^{\rho(x)}S(C_x,uv^{-1}) = \left(\frac{u}{v} - 1 \right)^{\rho(x)}
\sum_{m \in C_x \cap M} \left(\frac{u}{v}\right)^{{\rm deg}\, m}$$
on the right hand side of \ref{recur}.
Now we choose a lattice point $m \in C \cap M$,
collect terms containing $(u/v)^{{\rm deg}\, m}$
in right and left hand sides,
and use the equality (\ref{inv})
$$
\sum_{x(m) \leq x \leq \hat{1}}\left(\frac{u}{v} -1\right)^{\rho(x)}
G([x, \hat{1}],uv^{-1}) =
\left(\frac{u}{v} -1\right)^{\rho(x(m))}
\left(\frac{u}{v}\right)^{d-\rho(x(m))}
G([x(m),\hat{1}],u^{-1}v)
$$
on the right hand side.
By the duality  (\ref{duality0})
$$B([x(m),x]^*;u,v) = (-u)^{\rho(x) - \rho(m(x))}
B([x(m),x];u^{-1},v),$$
it remains to estabish  the recursive relation:
$$
\frac{(v-u)^{\rho(x(m))}}{uv}
\sum_{x(m) \leq x \leq \hat{1}} (-1)^{\rho(x)}(-u)^{\rho(x) - \rho(m(x))}
B([x(m),x];u^{-1},v) G([x,\hat{1}],uv)
 =
$$
$$
= \left(\frac{u}{v} -1\right)^{\rho(x(m))}
\frac{v^{d-1}}{u}\left(\frac{u}{v}\right)^{d-\rho(x(m))} G([x(m),\hat{1}],
u^{-1}v)
$$
which is equivalent to the recursive relation in
\ref{Q} after the substitution
$u^{-1}$ instead of $u$.
\hfill $\Box$

\section{Mirror duality}

Let $\overline{M}$ and $\overline{N} = {\rm Hom}(\overline{N}, {\bf Z})$
be dual to each other free abelian groups of rank $\overline{d}$,
$\overline{M}_{\bf R}$ and $\overline{N}_{\bf R}$
the real scalar extensions of $\overline{M}$ and $\overline{N}$,
$\langle *, *\rangle\;: \;
\overline{M} \times \overline{N} \rightarrow {\bf Z}$ the natural
pairing.

\begin{dfn}
{\rm \cite{batyrev-borisov1}
Let $C \subset \overline{M}_{\bf R}$ be a $\overline{d}$-dimensional
Gorenstein
cone. The cone $C$ is called {\em reflexive} if the dual cone
$C^* \subset N_{\bf R}$ is also Gorenstein; i.e., there exists
a lattice element $m_{C^*} \in M$ such that all vertices of the
supporting polyhedron $\Delta(C^*) = \{ z \in C^*\;:\;
\langle m_{C^*}, z \rangle =1 \}$ are contained in $M$. In this
case, we call $r = \langle m_{C^*}, n_{C} \rangle$ the {\em index} of $C$. }
\end{dfn}

\begin{dfn} {\rm
\cite{bat.dual} Let $M$ be a free abelian group of rank $d$.
A $d$-dimensional polyhedron in $M_{\bf R}$ with vertices in $M$ is
called {\em reflexive} if it can be identified with a
supporting polyhedron of some $(d+1)$-dimensional reflexive
Gorenstein cone of index $1$.}
\end{dfn}

Recall the definition of string-theoretic Hodge numbers of
an algebraic variety $X$ with at most Gorenstein toroidal
singularities
\cite{batyrev.dais}:

\begin{dfn}
{\rm \cite{batyrev.dais}
 Let $X = \bigcup_{i \in I} X_i$ be a $k$-dimensional
 stratified algebraic variety over ${\bf C}$ with at most Gorenstein
toroidal singularities such that for any $ i \in I$ the
singularities of $X$ along the stratum $X_i$ of codimension
$k_i$ are defined by a $k_i$-dimensional
finite rational polyhedral cone $\sigma_i$; i.e.,  $X$
is locally isomorphic to
$${\bf C}^{k-k_i} \times U_{\sigma_i}$$
at each  point $x \in X_i$ where $U_{\sigma_i}$ is a
$k_i$-dimensional  affine toric variety which is associated
with the cone $\sigma_i$ (see \cite{danilov}).
Then the polynomial
$$
E_{\rm st}(X;u,v) := \sum_{i \in I} E(X_i;u,v) \cdot
S(\sigma_i,uv)
$$
is called the {\em string-theoretic E-polynomial of $X$.} If we write
$E_{\rm st}(X; u,v)$ in form
$$
E_{\rm st}(X;u,v) = \sum_{p,q} a_{p,q} u^{p}v^{q},
$$
then the numbers $h^{p,q}_{\rm st}(X) := (-1)^{p+q}a_{p,q}$
are called the {\em string-theoretic Hodge numbers of $X$.}
}
\label{st-numbers}
\end{dfn}

\begin{rem}
{\rm Comparing with \ref{intersection} and \ref{e-form}, the
definition of the string-theoretic Hodge numbers looks as if there
were a complex ${ST}^{\bullet}_X$ whose hypercohomology groups
have natural Hodge structure which assumed to be pure if
$X$ is compact. We remark that the construction of such a  complex
${ST}^{\bullet}_X$ (an analog of the intersection complex) is
still an open problem.}
\end{rem}

Let $V = D_1 \cap \cdots \cap D_r$ be a generic
Calabi-Yau complete intersection of $r$ semi-ample divisors
$D_1, \ldots, D_r$ in a $d$-dimensional Gorenstein toric Fano variety
${\bf X}$ $(k \geq r)$. According to \cite{batyrev-borisov1},
there exists a $d$-dimensional reflexive polyhedron  $\Delta$ and
its decompostion
into a Minkowski sum
$$
\Delta = \Delta_1 +
\cdots + \Delta_r,$$
where each lattice polyhedron  $\Delta_i$  is the supporting
polyhedron for global sections of a semi-ample invertible
sheaf ${\cal L}_i \cong {\cal O}_{\bf X}(D_i)$
($i =1, \ldots, r$).

\begin{dfn}
{\rm \cite{borisov}
Denote by $E_1, \ldots, E_k$ the closures of $(d-1)$-dimensional
torus orbits in ${\bf X}$ and set $I := \{ 1, \ldots, k \}$.
A decompostion into a Minkowski sum
$\Delta = \Delta_1 + \cdots + \Delta_r$ as above is called a {\em
nef-partition} if there exists a decomposition of $I$ into a disjoint union
of $r$ subsets $I_j \subset I$ $(j =1, \ldots, r)$ such that
$$ {\cal O}(D_j) \cong  {\cal O}(\sum_{l \in I_j} E_l), \;\;
(j =1, \ldots, r)$$.}
\label{nef}
\end{dfn}

Now we put $\overline{M} = {\bf Z}^r \oplus M$, $\overline{d} = d+r$, and
define the $\overline{d}$-dimensional cone $C \subset \overline{M}_{\bf R}$
as
$$ C:= \{ (\lambda_1, \ldots, \lambda_r,
\lambda_1 z_1 + \cdots + \lambda_r z_r) \in \overline{M}_{\bf R}\; :\;
\lambda_i \in {\bf R}_{\geq 0}, \; z_i \in \Delta_i, \;
i =1, \ldots, r \}. $$
We extend the pairing $\langle \cdot, \cdot \rangle \; : M \times N
\rightarrow {\bf Z}$ to the pairing between $\overline{M}$ and
$\overline{N} := {\bf Z}^r \oplus N$ by the formula
$$\langle (a_1, \ldots, a_r, m), (b_1, \ldots, b_r,n) \rangle =
\sum_{i =1}^r a_i b_i + \langle m, n \rangle.$$

\begin{theo} {\rm \cite{borisov,batyrev-borisov1}} Let $\Delta =
\Delta_1 + \cdots + \Delta_r $ be a  nef-partition.
Then it defines canonically  a $d$-dimensional reflexive
polyhedron $\nabla \subset N_{\bf R}$
and a nef-partition $\nabla = \nabla_1 + \cdots + \nabla_r$
which are  uniquely determined by the property that
$$ C^*:= \{ (\lambda_1, \ldots, \lambda_r,
\lambda_1 z_1 + \cdots + \lambda_r z_r) \in \overline{N}_{\bf R}\; :\;
\lambda_i \in {\bf R}_{\geq 0}, \; z_i \in \nabla_i, \;
i =1, \ldots, r \}  $$
is the dual reflexive Gorenstein cone $C^*
\subset \overline{N}_{\bf R}$.
\label{nef-partition}
\end{theo}

\begin{dfn}
{\rm \cite{borisov} The nef-partition
$\nabla = \nabla_1 + \cdots + \nabla_r$ as in
\ref{nef-partition} is called {\em the
dual nef-partition}. }
\end{dfn}

We set
$${\bf Y}:= {\bf P}({\cal L}_1\oplus \cdots \oplus
{\cal L}_r).
$$
Recall  the standard construction of the reduction of complete
intersection $V \subset {\bf X}$
to a hypersurface $\tilde{V} \subset {\bf Y}$ \cite{batyrev-borisov1}.
Let $\pi$ be the canonical projection ${\bf Y} \rightarrow {\bf X}$ and
${\cal O}_{\bf Y}(-1)$ the tautological Grothendieck sheaf on ${\bf Y}$. Since
$$\pi_*{\cal O}_{\bf Y}(1) = {\cal L}_1\oplus \cdots \oplus
{\cal L}_r,$$
we obtain the isomorphism
$$ H^0({\bf Y}, {\cal O}_{\bf Y}(1)) \cong H^0({\bf X}, {\cal L}_1) \oplus
\cdots \oplus H^0({\bf X}, {\cal L}_r).$$
Assume that $D_i$ is the set of zeros of a global section $s_i \in
H^0({\bf X}, {\cal L}_i)$ $(1 \leq i \leq r)$. We define
$\tilde{V}$ as the zero set of the global section
$s \in H^0({\bf Y}, {\cal O}_{\bf Y})$ which corresponds to the
$r$-tuple $(s_1, \ldots, s_r)$ under above isomorphism. Our main interest
is the following standard property (\cite{batyrev-borisov1}):

\begin{prop}
The restriction of $\pi$ on ${\bf Y} \setminus \tilde{V}$ is a locally trivial
${\bf C}^{r-1}$-bundle in Zariski topology over ${\bf X} \setminus V$.
\label{bundle}
\end{prop}

Let us set
$${\bf P} = {\rm Proj}\, \bigoplus_{i \geq 0} H^0({\bf Y},
{\cal O}_{\bf Y}(i)).$$
The following statement is contained in \cite{batyrev-borisov1}:

\begin{prop}
The tautological sheaf ${\cal O}_{\bf Y}(1)$ is semi-ample and the natural
toric morphism
$$\alpha \;: \; {\bf Y} \rightarrow {\bf P}$$
is crepant. Moreover, ${\cal O}_{\bf Y}(r)$ is the anticanonical sheaf  of
${\bf Y}$, ${\bf P}$ is a Gorenstein toric Fano variety , and
$\overline{Z} : = \alpha(\tilde{V})$ is an ample hypersurface in ${\bf P}$.
\label{alpha}
\end{prop}

There is the following explicit formula for
$E_{\rm st}(V;u,v)$ in terms of
$E_{\rm st}({\bf P};u,v)$ and
$E_{\rm st}(\overline{Z};u,v)$:

\begin{theo}
$$E_{\rm st}(V;u,v)
=((uv-1)((uv)^r-1)^{-1})E_{\rm st}({\bf P};u,v)
 - (uv)^{1-r}E_{\rm st}({\bf P} \setminus \overline{Z};u,v).$$
\label{V-form}
\end{theo}

\noindent
{\em Proof. } Since $V$ is transversal to all toric
strata in ${\bf X}$ we have:
$$ E_{\rm st}(V;u,v)=
E_{\rm st}({\bf X};u,v) -  E_{\rm st}({\bf X} \setminus V;u,v).$$
Using the ${\bf CP}^{r-1}$-bundle structure of ${\bf Y}$ over
${\bf X}$, we obtain:
$$
E_{\rm st}({\bf X};u,v) = ((uv)^r-1)^{-1}(uv-1)E_{\rm st}({\bf
Y};u,v).
$$
By \ref{bundle}, we also have
$$ E_{st}({\bf X} \setminus V;u,v)
= (uv)^{1-r} E_{\rm st}({\bf Y} \setminus \tilde{V};u,v).$$
Since birational crepant toric morphisms do not change string-theoretic
Hodge numbers (see \cite{batyrev.dais}), by \ref{alpha}, we conclude
$$ E_{\rm st}({\bf Y};u,v)= E_{\rm st}({\bf
P};u,v), \;\;
E_{\rm st}({\bf Y} \setminus \tilde{V};u,v) =
E_{\rm st}({\bf P} \setminus \overline{Z};u,v) .$$
\hfill $\Box$

\begin{dfn}
{\rm Let $C \subset \overline{M}_{\bf R}$ be a reflexive Gorenstein cone,
$C^* \subset {\overline{N}}_{\bf R}$ the dual reflexive
Gorenstein cone. We define
$$
\Lambda(C,C^*):= \{ (m,n) \in \overline{M} \oplus
\overline{N}\;: \; m \in C, \;
n \in C^*,\;\;\mbox{\rm and}\;   \; \langle m, n \rangle = 0 \}.$$
}
\end{dfn}

\begin{dfn}
{\rm Let $(m,n)$ be an element of $\Lambda(C,C^*)$. We define
the Eulerian poset $P_{(m,n)}$ as the subset of all faces $C_x
\subset C$ such that $C_x$ contains $m$ and
$\langle z,n \rangle = 0$ for all $z \in C_x$. We denote by
$\rho(x^*(n))$ the dimension of
the intersection of $C$ with the hyperplane $\langle z,n \rangle = 0$.
}
\label{n-rho}
\end{dfn}

\begin{rem}
{\rm The dual Eulerian poset $P^*_{(m,n)}$ can be identified with the
subset of all faces $C_x^* \subset C^*$ such that $C_x^*$ contains $n$ and
$\langle m,z \rangle = 0$ for all $z \in C_x^*$. }
\end{rem}

\begin{theo}
Let us set $\overline{d} = d + r$ and
$$A_{(m,n)}(u,v) =
\frac{(-1)^{\rho(x^*(n))}}{(uv)^r}
(v-u)^{\rho(x(m))}B(P_{(m,n)}^*;u,v)(uv -1)^{\overline{d}-\rho(x^*(n))}.
$$
Then
$$
E_{\rm st}(V; u,v)=
\sum_{(m,n) \in \Lambda(C,C^*)}
\left(\frac{u}{v}\right)^{{\rm deg}\,m} A_{(m,n)}(u,v)
\left(\frac{1}{uv}\right)^{{\rm deg}\,n}
$$
\label{st.formula}
\end{theo}

\noindent
{\em Proof.}
By Definition \ref{st-numbers},
$$
E_{\rm st}({\bf P},u,v) = \sum_{\hat{0} < x \leq \hat{1}}
(uv - 1)^{\rho(x) -1} S(C_x^*,uv) $$
$$ =
\sum_{\hat{0} < x \leq \hat{1}}
(uv - 1)^{\rho(x) -1} (uv -1)^{\overline{d} - \rho(x)} T(C_x^*,(uv)^{-1})
$$
$$
= (uv -1)^{\overline{d} -1} \sum_{\hat{0} < x \leq \hat{1}}
  \left(\sum_{n \in
Int(C_x^*) \cap {\overline{N}}} (uv)^{-{\rm deg}\,n} \right) =
(uv -1)^{\overline{d}-1} \sum_{n \in
\partial C^* \cap {\overline{N}}} (uv)^{-{\rm deg}\,n},
$$
where $\partial C^* = C^* \setminus Int(C^*)$ is the boundary of $C^*$.
Since ${\overline{N}} \cap Int(C^*) = p+ {\overline{N}}\cap C^*$ and ${\rm
deg}\, p = r$, we conclude:
$$
E_{\rm st}({\bf P},u,v) =
(1-(uv)^{-r})(uv -1)^{\overline{d}-1} \sum_{n \in  C^* \cap {\overline{N}}}
(uv)^{-{\rm deg}\,n}$$
$$ =
((uv)^{r}-1)(uv -1)^{\overline{d}-1} \sum_{n \in Int( C^*) \cap {\overline{N}}}
(uv)^{-{\rm deg}\,n}.
$$
On the other hand,
$$
E_{\rm st}({\bf P} \setminus \overline{Z};u,v) =
E_{\rm st}({\bf P};u,v) - E_{\rm st}(\overline{Z};u,v). $$
By Definition \ref{st-numbers}  and Theorem \ref{f-formula},
$$
E_{\rm st}(\overline{Z};u,v) =
\sum_{\hat{0} < x \leq \hat{1}}
\left( \frac{(uv-1)^{\rho(x)-1}}{uv} \right) S(C_x^*,uv)$$
$$+ \sum_{\hat{0} < x \leq \hat{1}}
\left( \frac{(-1)^{\rho(x)}}{uv} \sum_{m \in C_x \cap {\overline{M}}}
(v-u)^{\rho(x(m))}B([x(m),x]^*;u,v)
\left(\frac{u}{v}\right)^{{\rm deg}\,m}\right)
S(C_x^*,uv)$$
$
= (uv)^{-1} E_{\rm st}({\bf P};u,v) + $
$$ +
\sum_{\hat{0} < x \leq \hat{1}}
\left( \frac{(-1)^{\rho(x)}}{uv} \sum_{m \in C_x \cap {\overline{M}}}
(v-u)^{\rho(x(m))}B([x(m),x]^*;u,v) \left(\frac{u}{v}\right)^{{\rm deg}\,m}
\right) S(C_x^*,uv).
$$
By \ref{V-form},
$$E_{\rm st}(V;u,v) =
((uv-1)((uv)^r-1)^{-1} - (uv)^{1-r} + (uv)^{-r})
E_{\rm st}({\bf P},u,v) $$
$$ +  \sum_{\hat{0} < x \leq \hat{1}}
\left( \frac{(-1)^{\rho(x)}}{(uv)^r} \sum_{m \in C_x \cap {\overline{M}}}
(v-u)^{\rho(x(m))}B([x(m),x]^*;u,v) \left(\frac{u}{v}\right)^{{\rm deg}\,m}
\right) S(C_x^*,uv)
$$
$$
= (uv)^{-r}(uv -1)^{\overline{d}} \sum_{n \in  Int(C^*)
\cap {\overline{N}}} (uv)^{-{\rm deg}\,n}
$$
$$
+  \sum_{\hat{0} < x \leq \hat{1}}
\left( \frac{(-1)^{\rho(x)}}{(uv)^r} \sum_{m \in C_x \cap {\overline{M}}}
(v-u)^{\rho(x(m))}B([x(m),x]^*;u,v) \left(\frac{u}{v}\right)^{{\rm deg}\,m}
\right) S(C_x^*,uv)
$$
$$
=  \sum_{\hat{0} \leq  x \leq \hat{1}}
\left( \frac{(-1)^{\rho(x)}}{(uv)^r} \sum_{m \in C_x \cap {\overline{M}}}
(v-u)^{\rho(x(m))}B([x(m),x]^*;u,v)
\left(\frac{u}{v}\right)^{{\rm deg}\,m}
\right) S(C_x^*,uv).
$$
It remains to use  the formula
$$
S(C_x^*,uv)= (uv-1)^{\overline{d} -\rho(x)} \sum_{n \in
Int(C_x^*)\cap {\overline{N}}} (uv)^{-{\rm
deg}\,n} \;\;\; ( \hat{0} \leq  x \leq \hat{1})
$$
and notice that
$ \rho(x)=  \rho(x^*(n))$ if  $n$ is an interior lattice
point of $C_x^*$ (see  \ref{n-rho}).
\hfill
$\Box$

\begin{theo}
Let $V$ be a $(d-r)$-dimensional  Calabi-Yau complete
intersection  defined by a nef-partition
$\Delta = \Delta_1 + \cdots + \Delta_r$,
$W$ a $(d-r)$-dimensional Calabi-Yau complete intersection defined by the
dual nef-partition $\nabla = \nabla_1 + \cdots + \nabla_r$. Then
$$E_{\rm st}(V;u,v) = (-u)^{d-r}E_{\rm st}(W;u^{-1},v),$$
i.e.,
$$
h^{p,q}_{\rm st}(V) = h^{d-r-p,q}_{\rm st}(W)\;\; 0 \leq p, q
\leq d-r.
$$
\label{duality1}
\end{theo}

\noindent
{\em Proof.} If we use the duality between two
$\overline{d}$-dimensional reflexive Gorenstein cones
$C \subset \overline{M}_{\bf
R}$ and $C^* \subset \overline{N}_{\bf R}$ \ref{nef-partition}, then the
statement of Theorem follows immediatelly from the explicit
formula in \ref{st.formula} and from the duality for
$B$-polynomials \ref{duality0}.

\end{document}